\documentclass[12pt,oneside]{amsart}

\usepackage{latexsym,amssymb,amsmath,graphicx,tkz-graph,hyperref,verbatim}
\usepackage{subfig}
\usepackage{caption}
\numberwithin{equation}{section}
\newtheorem{theorem}{Theorem}[section]
\newtheorem{lemma}[theorem]{Lemma}

\newtheorem{question}[theorem]{Question}
\theoremstyle{definition}
\newtheorem{definition}[theorem]{Definition}

\newtheorem{example}[theorem]{Example}

\newtheorem{algorithm}[theorem]{Algorithm}

\DeclareMathOperator{\lcm}{lcm}
\DeclareMathOperator{\Sym}{Sym}

\newenvironment{steps}
{\begin{enumerate}
	}
{\end{enumerate}}

\begin{document}


\title{Revisiting the spreading and covering numbers}
\thanks{Last updated: \today}

\author{Ben Babcock}
\author{Adam Van Tuyl}
\address{Department of Mathematical Sciences\\
Lakehead University, Thunder Bay, ON, P7B 5E1, Canada}
\email{bababcoc@lakeheadu.ca, avantuyl@lakeheadu.ca}


\begin{abstract}
We revisit the problem of computing the spreading and covering numbers.
We 
show a connection between some of the spreading
numbers and the number of non-negative integer $2 \times 2$ matrices
whose entries sum to $d$, and we construct an algorithm to compute
improved upper bounds for the covering numbers.
\end{abstract}

\maketitle


\section{Introduction}\label{sec:intro}

Let $R = k[x_1,\ldots,x_n]$ be a polynomial ring over a field $k$.  For any
non-negative integer $d$, let $M_d$ be the set of all monomials of degree
$d$ in $R$.  For any subset $W \subseteq M_d$, let 
\[R_1W = \{x_im ~:~ m \in W ~~\text{and}~~ 1 \leq i \leq n \}. \]   
For any $W \subseteq M_d$, we always have
$|R_1W| \leq n|W|$ and $R_1W \subseteq M_{d+1}$.  

We are interested
in finding subsets $W$ where either $|R_1W| = n|W|$ or $R_1W = M_{d+1}$.  
We define the {\it spreading number} to be
\[\alpha_n(d) = \max\{ ~|W| ~:~ W \subseteq M_d ~~~\text{and}~~~ |R_1W| = n|W|\}.\]
The terminology is derived
from the fact that the elements of $R_1W$ are ``spread'' out
in $M_{d+1}$.  Similarly, the {\it covering number} is defined to be
\[\rho_n(d+1) = \min\left\{ ~|W| ~:~ W 
\subseteq M_d ~~\text{and}~~~ R_1W = M_{d+1} \right\}.\]
In this case the elements of $R_1W$ ``cover'' the elements of $M_{d+1}$.

Geramita, Gregory, and Roberts introduced
$\alpha_n(d)$ and $\rho_n(d+1)$
to study the Ideal Generation Conjecture for a set of generic points in
$\mathbb{P}^n$ (see \cite[Theorem 4.7]{Geramita86}).
When $n=1$, it is trivial to show that $\alpha_1(d)=\rho_n(d+1) = 1$ for
all $d$.  When $n=2$, $\alpha_2(d) = \lfloor \frac{d}{2} \rfloor +1$ and
$\rho_2(d+1) = \lceil\frac{d}{2} \rceil +1$. 
Geramita, et al. gave
exact values for $\alpha_n(d)$ for all $d$ when $n=3$ or $4$,
and some scattered results and bounds on other values.  
Curtis \cite{Curtis95} later found a formula
for $\rho_3(d)$ for all $d$ and an improved lower bound on $\rho_4(d)$.
Using techniques from linear programming, Hulett and Will \cite{Hulett99}
improved these lower bounds on $\rho_4(d)$.  Carlini,
H\`a, and the second author \cite{Carlini01} later reformulated the
problem by constructing simplicial complexes whose dimensions
were related to either $\alpha_n(d)$ or $\rho_n(d+1)$. 

Surprisingly, computing new exact values of $\alpha_n(d)$ and $\rho_n(d+1)$
remains elusive.  However, we present two new contributions: 
1)  a new connection between the numbers
$\alpha_4(d)$ and 
the number of integer matrices with a specific property; and 2)
a new greedy algorithm which gives upper bounds on $\rho_n(d+1)$ 
that improves upon known bounds.  Hopefully these observations will
be of use for future attacks on computing $\alpha_n(d)$ and $\rho_n(d+1)$.


\section{Preliminaries}\label{sec:prelim}

We begin by translating our problem of computing $\alpha_n(d)$
and $\rho_n(d+1)$ into a graph theory problem.
Fix positive integers $n$ and $d$.  Let $S_n(d)$ denote the
graph whose vertex set is the set of monomials $M_d$ in $R = k[x_1,\ldots,x_n]$,
and two vertices $m_i, m_j$ are adjacent if and only if
$\deg\big(\lcm(m_i, m_j)\big) = d + 1.$  
We abuse notation and use $M_d$ to denote both the vertices of $S_n(d)$ and the
set of monomials of degree $d$ in $R = k[x_1,\ldots,x_n]$.
We denote the number of vertices of $S_n(d)$ by $v_d(n)$;  it is clear
that $v_d(n) = \binom{n+d-1}{d}$.

\begin{definition}
A subset $V \subseteq M_d$ is an {\it independent set} if any two distinct 
elements of $V$ are not adjacent; $V$ is a {\it maximal 
independent set} if it is not properly contained in any larger independent set.
\end{definition}

\begin{definition} A subset of $M_d$ in which any two vertices are 
adjacent is called a {\it clique}. 
If $C_1,\ldots,C_t$ are cliques,
we say they form a {\it clique cover} of $S_n(d)$ if 
$C_1 \cup \cdots \cup C_t = M_d$.
For any monomial $m$ of degree $d - 1$, an 
{\it upward clique} is the clique consisting of the vertices 
$mx_i \in M_d$ for all $x_i \in \{x_1,\ldots,x_n\}$.
\end{definition}

As shown in \cite{Geramita86}, $\alpha_n(d)$ and $\rho_n(d+1)$
are equivalent to an invariant of  $S_n(d)$:

\begin{lemma} \label{graphtranslate}
 With the notation as above
\begin{enumerate}
\item[$(i)$] $\alpha_n(d)$ is the cardinality of the largest maximal independent set of $S_n(d)$.
\item[$(ii)$] $\rho_n(d+1)$ is the minimum cardinality of an upward clique cover of 
the vertices of $S_n(d + 1)$.
\end{enumerate}
\end{lemma}

\begin{example}
If we consider
$S_3(3)$ (see Figure 1), then 
$x_1^2x_2$ and $x_1^2x_3$ are adjacent but $x_1^2x_2$ and $x_2^2x_3$ are not.
\begin{figure}\label{graph1}
\begin{tikzpicture}[node distance=1cm]
	\GraphInit[vstyle=Classic]
	\SetVertexMath

	\Vertex[Lpos=90]{x_1^3} \SOWE[Lpos=180](x_1^3){x_1^2x_2} \SOEA(x_1^3){x_1^2x_3}
	\Edges(x_1^3,x_1^2x_2,x_1^2x_3,x_1^3)

	\SOWE[Lpos=180](x_1^2x_2){x_1x_2^2} \SOEA[Lpos=90](x_1^2x_2){x_1x_2x_3} \SOEA(x_1^2x_3){x_1x_3^2}
	\Edges(x_1^2x_2, x_1x_2^2, x_1x_2x_3, x_1x_3^2, x_1^2x_3, x_1x_2x_3, x_1^2x_2)

	\SOWE[Lpos=210](x_1x_2^2){x_2^3} \SOEA[Lpos=270](x_1x_2^2){x_2^2x_3} \SOEA[Lpos=270](x_1x_2x_3){x_2x_3^2} \SOEA[Lpos=300](x_1x_3^2){x_3^3}
	\Edges(x_1x_2^2, x_2^3, x_2^2x_3, x_2x_3^2, x_3^3, x_1x_3^2, x_2x_3^2, x_1x_2x_3, x_2^2x_3, x_1x_2^2)
\end{tikzpicture}
\caption{The graph $S_3(3)$}
\end{figure}
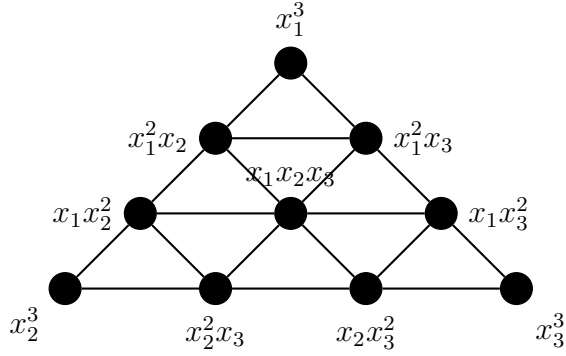
The graph $S_3(3)$ has $\alpha_3(3) = 4$ because
$\{x_1^3,x_2^3,x_3^3,x_1x_2x_3\}$ forms a maximal independent set.  
Also $\rho_3(2+1) = 4$ because $C_1 = \{x_1^3,x_1^2x_2,x_1^2x_3\}$,
$C_2 = \{x_2^3,x_1x_2^2,x_3x_2^2\}$, $C_3 = \{x_3^3,x_1x_3^2,x_2x_3^2\}$,
and $C_4 = \{x_1^2x_2,x_1x_2^2,x_1x_2x_3\}$ form a minimal upward clique cover.
\end{example}

Computing the size of a largest 
maximal independent set or a minimum clique cover 
of a graph are both NP-hard problems.  This explains,
in part, why
computing $\alpha_n(d)$ and $\rho_n(d+1)$ is so difficult.


\section{A053307}\label{sec:sequence}

We demonstrate a relation between the sequence $\alpha_4(d)$ 
and a known integer sequence which is denoted in the OEIS as 
A053307 \cite{OEIS}.
In \cite{Geramita86} we find the explicit formula:
\[
\alpha_4(d) = \left\{
\begin{tabular}{ll}
$\frac{v_4(d)}{4}$ & \mbox{for $d$ odd} \\
$\frac{v_4(d)}{4} + \frac{3d + 6}{8}$ & \mbox{for $d$ even.}
\end{tabular}
\right.
\] 

\begin{theorem}
For all $d \geq 0$, $\alpha_4(d)$ 
equals the number of non-negative integer $2\times2$ 
matrices with sum of entries equal to $d$, under row and column permutations.
\end{theorem}

\begin{proof}
Recalling that $v_4(d) = \binom{d + 3}{3}$, it follows that
\begin{align*}
\alpha_4(2d+1) 	&= 
\frac{\binom{2d + 4}{3}}{4} = \frac{(2d + 4)(2d + 3)(2d + 2)}{3!(4)} = \frac{(d+2)(2d+3)(d+1)}{6}.
\end{align*}
Similarly, 
\begin{align*}
\alpha_4(2d) 	&= \frac{\binom{2d+3}{3}}{4} + \frac{6d + 6}{8} = \frac{8d^3 + 24d^2 + 40d+ 24}{24} = \frac{(d+1)^3 + 2(d+1)}{3}.
	\end{align*}
The OEIS reveals that $\alpha_4(2d+1) =  \text{A000330}(d + 1)$;
the sequence \text{A000330} is the sequence whose
$i$-th term is given by $0^2+1^2+\cdots+i^2$.
Similarly, $\alpha_4(2d) = \text{A006527}(d+1)$, 
the sequence whose $i$-th term is $(i^3+2i)/3$.
So, $\alpha_4(d)$ is an interleaved sequence.

Let $a(d)$ be the number of non-negative integer $2\times2$ matrices 
with sum of 
entries equal to $d$, under row and column permutations.  
The OEIS lists this sequence as A053307, and contains
a comment,  attributed to Paul Barry, that the 
integer sequence \text{A053307} is also the interleaved sequence
of  \text{A000330} and \text{A006527}, i.e., 
$\text{A053307}(2d+1) =  \text{A000330}(d + 1)$ and 
$\text{A053307}(2d) =  \text{A006527}(d+1)$.   The conclusion follows from
this observation.

Since no proof is given for Barry's comment, 
we sketch out why this is indeed
the case.  The generating function for $\text{A053307}$ is listed in the OEIS 
as
\[\frac{t^2-t+1}{(1-t^2)^2(1-t)^2}.\]
Multiplying the top and bottom of this expression by $(1+t)^2$ gives
\[ \frac{(t^2-t+1)(1+t)^2}{(1-t^2)^2(1-t)^2(1+t)^2} =
\frac{t^4+1}{(1-t^2)^4} + \frac{t(t^2+1)}{(1-t^2)^4}.\]
It follows that $\text{A053307}(2d)$ equals the coefficient of $t^{2d}$
of $\frac{t^4+1}{(1-t^2)^4}$, and $\text{A053307}(2d+1)$ equals the
coefficient of $t^{2d+1}$ in the other rational function.
Now the rational  function $\frac{(t^2+1)}{(1-t)^4}$ is the generating function
of A006527 (this is slightly different than what is listed in the
OEIS because we want the sequence to start with  $1$, not $0$, 
so we have dropped
the extra multiple $t$).  Replacing $t$ with $t^2$ gives the 
first rational function on the right hand side, which means
$\text{A006527}(d+1) = \text{A053307}(2d)$.  A similar analysis
using $\frac{t(t+1)}{(1-t)^4}$, the generating
function of A000330, will complete the proof.
\end{proof}

Even though the sequence A053307 and $\alpha_4(d)$ are related, it is not
immediately apparent why they are linked, thus suggesting the following question:

\begin{question}
Is there an explicit bijection between the the maximal independent
sets of $S_4(d)$ and the number of non-negative integer $2\times2$ 
matrices with sum of entries equal to $d$, under row and column permutations?
\end{question}

The correspondence may be a result of the two interleaved 
sequences that make up A053307.
Explaining the relationship between 
$\alpha_4(d)$ and A053007 could open up new techniques 
for computing the spreading and covering numbers.


\section{A Greedy Algorithm for bounding $\rho_n(d)$}\label{sec:greedy-bounds}

We use the symmetry of the graph $S_n(d)$
to describe a greedy algorithm that bounds from above $\rho_n(d)$.  We give
evidence
that our algorithm improves on known bounds.

\subsection{The Algorithm}
By Lemma \ref{graphtranslate} $\rho_n(d)$ is the cardinality of the minimum 
upward clique cover of $S_n(d)$.  We give a greedy algorithm
that constructs an upward clique cover.  Roughly speaking,
at each step, the algorithm picks an upward clique for any
vertex that has not been covered.  The number of upward cliques in our
cover forms our bound on $\rho_n(d)$. 

We begin with some observations.  By definition, every upward clique is 
uniquely identified with a monomial from $M_{d-1}$. For vertices of $S_n(d)$ 
that consist of more than one 
indeterminate, many factorizations into a degree $d-1$ monomial and a variable
are possible; e.g., $x_1x_2^2$ can be 
written as $(x_1)x_2^2$ or $(x_1x_2)x_2$.  However, for monomials 
of the form $x_i^{d}$,  there is one such factorization,
that is $(x_i^{d-1}){x_i}$.  Thus, $x_i^{d}$ belongs only to the upward 
clique identified with $x_i^{d-1}$. These unique upward 
cliques containing each $x_i^{d}$ must therefore be in any clique cover of 
$S_n(d)$, so we can use them as the our initial set.

Aside from our choice of initial members of the cover, we wish to take
into account the symmetry of  $S_n(d)$. 
Let $\Sym(n)$ denote the symmetric group on the set 
$\{1,2,\ldots,n\}$.
For any $\mathbf{x^a} = x_1^{a_1}\cdots{x}_n^{a_n} \in M_d$ and 
$\sigma\in \Sym(n)$, let $\sigma(\mathbf{x^a})$ be
the monomial obtained by permuting the indices $1,\ldots,n$ according to the permutation
$\sigma$. 
This operation preserves many properties of sets of 
vertices; e.g., independent sets and clique covers are both unaffected.

We use $\Sym(n)$ to create  orbits of the vertices of $S_n(d)$; that is,
for any $m \in M_d$, the {\it orbit of $m$} is the set
$\{\sigma(m) ~|~ \sigma \in \Sym(n)\}.$
Since elements 
of $\Sym(n)$ do not alter the exponents of a monomial, only the order of the 
exponents relative to the indeterminates, the orbit of $m$ is also the set of
 all permutations of the exponents of $m$. By definition, the exponents of any
 $m \in M_d$ always sum to $d$, and therefore the orbits of $S_n(d)$
are in an one-to-one correspondence with 
the integer partitions of $d$ of length at most $n$.  
We can write orbits 
as vectors in $\mathbb{N}^n$, and in this form it is easy to 
determine whether an
orbit is an independent set, a clique, or neither by examining 
the entries in the vector. 
We will order our list of orbits with respect to the
reverse lexicographical order, that is, if $\alpha,\beta \in \mathbb{N}^n$,
then $\alpha \geq_{rlex} \beta$ if the last non-zero entry of $\alpha-\beta$
is negative.
Iterating over the list of orbits of 
$S_n(d)$ in reverse lexicographical order will help us in computing
an upper bound on $\rho_n(d)$.

We now present our algorithm that returns a minimal upward clique cover;
$\rho_n(d)$ is bounded above by the number of cliques in this cover.
\medskip

\hrule

\begin{algorithm}\label{alg:covering-orbits}
Compute an upper bound for $\rho_n(d)$.

\noindent
\textbf{Input:} $n,d$ --- The number of variables  and degree of
 monomials, respectively. \\
\textbf{Output:} A minimal upward clique cover of $S_n(d).$
\smallskip

\hrule
\smallskip

\begin{steps}
	\item Initialize our cover $\mathcal{C}$ with  the set of 
upward cliques that contain $x_i^d$.
	\item Obtain a list, $L$, of the orbits of $S_n(d)$, where each orbit is 
represented as a vector in $\mathbb{N}^n$. Sort the list in reverse lexicographical order.
	\item Iterate over $L$. For each orbit $O \in L$, iterate over the 
vertices $v \in O$. If $v$ is covered, continue. If not, iterate over the 
upward cliques containing $v$. Select the upward clique that contains the fewest number of vertices already covered, and add it to $\mathcal{C}$.
	\item \label{lst:compute-freq} For each $v \in M_d$, 
compute its frequency, i.e., the number of upward cliques that contain it, in $\mathcal{C}$.
	\item 
 Iterate over the elements of the $\mathcal{C}$.  If an upward clique does not 
contain a vertex of frequency 1---all its vertices are represented by other 
cliques as well---then it is not essential to the cover, so discard it.  
Repeat this step until we complete an iteration without discarding any cliques.
	\item Return $\mathcal{C}$ as a minimal cover.
\end{steps}
\end{algorithm}
\hrule

\subsection{Comparison to Known Bounds}\label{sec:results}

We compare the known bounds for $\rho_n(d)$ to the output
of Algorithm \ref{alg:covering-orbits}.
Geramita, Gregory, and Roberts proved:

\begin{theorem}[{\cite[Theorem 5.2 and Proposition 5.9]{Geramita86}}]\label{thm:ggr-bounds}
	For all $n \geq 2, d \geq 2$, 
\[ \frac{v_n(d)}{n} \leq \alpha_n(d) \leq \rho_n(d) \leq \frac{v_n(d)}{n} + \frac{n - 1}{n}v_{n-1}(d) ~~\mbox{where $v_n(d) = \binom{n+d-1}{d}$}.\]
\end{theorem}

Hulett and Will improved the bounds on $\rho_4(d)$:

\begin{theorem}[{\cite[Theorems 4.1 and 4.2]{Hulett99}}]\label{thm:hw-upper}
	For all $d \geq 5$, 
	\begin{enumerate}
		\item[$(i)$] if $d$ is odd, $\rho_4(d) \leq (d^3 + 15d^2 - 61d + 261)/24$, or
		\item[$(ii)$] if $d$ is even, $\rho_4(d) \leq (d^3 + 15d^2 - 34d + 240)/24$.
	\end{enumerate}
\end{theorem}


We first consider the values of $\rho_4(d)$. In 
Table \ref{tab:covering-comparison}, {\bf GGR} refers to the 
upper bound for $\rho_4(d)$  in 
Theorem \ref{thm:ggr-bounds} and {\bf HW} refers to the bounds 
from Theorem \ref{thm:hw-upper}, while 
{\bf \ref{alg:covering-orbits}} refers to the bounds found using  
Algorithm \ref{alg:covering-orbits}.

\noindent
\begin{minipage}{.45\linewidth}
\centering

\begin{tabular}{||r|r|r|r||}
\hline
$d$ & 	 GGR & 	HW & \ref{alg:covering-orbits} \\
\hline
\hline
$5$ & 	30 & 	19 & 	19 \\
$6$ & 	42 & 	33 & 	29 \\
$7$ & 	57 & 	38 & 	40 \\
$8$ & 	75 & 	60 & 	55 \\
$9$ & 	97 & 	69 & 	74 \\
$10$ & 	 121 & 	 100 & 	96 \\
$11$ & 	 150 & 	 114 & 	 122 \\
$12$ & 	 182 & 	 155 & 	 147 \\
$13$ & 	 219 & 	 175 & 	 185 \\
$14$ & 	 260 & 	 227 & 	 223 \\
$15$ & 	 306 & 	 254 & 	 275 \\
\hline
\end{tabular}
\captionof{table}{Comparison of upper bounds for $\rho_4(d).$}
\label{tab:covering-comparison}

\noindent

\end{minipage}\hfill
\begin{minipage}{.45\linewidth}
\centering

\begin{tabular}{||r|r|r||}
\hline
$d$ & 	 GGR & \ref{alg:covering-orbits} \\  
\hline
\hline
$6$ & 110	& 61\\
$7$ & 162	& 94\\
$8$ & 231	& 142\\
$9$ & 319	& 209\\
$10$ & 429	& 285\\
$11$ & 565	& 392\\
$12$ & 728 &515 \\
$13$ & 924 &671 \\
$14$ & 1156&872\\
\hline
\end{tabular}
\captionof{table}{Comparison of upper bounds for $\rho_5(d).$}
\label{tab:covering-comparison-2}

\end{minipage}

The output of 
Algorithm \ref{alg:covering-orbits}  is quite close to HW. 
In fact, it seems that for even $d$ our bounds are equal or better, 
while the reverse is true for odd $d$. This pattern holds for at least 
$d \leq 24$, with the  exception of $d = 22$. 
We are not certain why this is the case. 

While the  HW bound holds only for $n = 4$, our algorithm works for all 
$n \geq 2$.   When tested against GGR for small values of $d$ for $n = 5,6$,
Algorithm \ref{alg:covering-orbits} consistently performs better.
Refer to Table \ref{tab:covering-comparison-2} for a comparison when $n=5$.
We hope  this provides a useful example of how one can use 
the structure and symmetry of $S_n(d)$ 
along with a greedy algorithm to improve bounds on $\rho_n(d)$.  It
also suggests that the bounds of GGR are far from optimal.

\subsection{Comments on implementation}
Some of the  computations were performed in \emph{Macaulay2} 1.3.1 \cite{Mt}
with 4 GB of 
memory allocated to 1 CPU and 1 node on SHARCNET's Saw 
cluster.\footnote{\url{https://www.sharcnet.ca/my/systems/show/41}}
Other computations ran on the Kraken 
cluster\footnote{\url{https://www.sharcnet.ca/my/systems/show/69}} 
and used up to 16 GB of memory in \emph{Macaulay2} 1.4. 
Readers interested in our code can visit our websites\footnote{\url{https://github.com/tachyondecay/spreading-covering-numbers/}
\newline
\url{http://flash.lakeheadu.ca/~avantuyl/research/SpreadCover_Babcock_VanTuyl.html}}.
The run times 
are taken from \emph{Macaulay2}'s \texttt{time} function.
The algorithm does not consume much memory, but as one might expect, 
as $d$ increases the 
computational time increases significantly. As a result, 
we found it difficult to compute bounds beyond $d > 10$.
When $n=4$, the largest $d$ for which we could compute a bound
was for $\rho_4(24)$.  In this case, the computation took 83051.40 seconds.


\section{An additional (unsuccessful) attack}

We end with a description of a theoretic approach
for bounding the numbers $\alpha_n(d)$
using commutative algebra.
While present computing power does not
enable us to apply this approach, we record this method for future attacks.

Given a finite simple graph $G$ with vertex set $V_G = \{z_1,\ldots,z_t\}$
and edge set $E_G$, 
the {\it edge ideal} of $G$ is 
$I(G) = (z_iz_j ~|~ \{z_i,z_j\} \in E_G) \subseteq T=k[z_1,\ldots,z_t].$
Some of the graph invariants of $G$ are encoded
into the algebraic invariants of $I(G)$.  For example,  it is known 
(e.g., see \cite{Vi}) that the Krull dimension of $T/I(G)$, denoted
$\dim T/I(G)$, equals $\alpha(G)$, the {\it independence number} of
$G$, that is, is the cardinality of the maximum independent set.
When $G = S_n(d)$, it follows by Lemma \ref{graphtranslate}
that $\alpha(S_n(d))= \alpha_n(d)$ and  thus $\alpha_n(d) 
= \dim T/I(S_n(d))$ where $T = k[z_m ~|~ m \in M_d]$.

To compute or bound $\alpha_n(d)$, it therefore suffices to compute
or bound the dimension of a ring.  One approach, therefore,
is to make use of following lemma:

\begin{lemma}  \label{phsop}
Let $L_1,\ldots,L_t$ be any $t$ linear forms of $T$.
Then
\[\dim T/(I(S_n(d)),L_1,\ldots,L_t) + t \geq \dim T/I(S_n(d)) = \alpha_n(d).\] 
\end{lemma}

\begin{proof}
This follows from the more
general fact that for any homogeneous ideal $I$ in $T$
and linear form $L \in T$,  then $\dim T/(I,L) \geq \dim T/I -1$.
\end{proof}

A strategy to bound $\alpha_n(d)$ is to find linear forms
$L_1,\ldots,L_t$ so that the computation of 
$\dim T/(I(S_n(d)),L_1,\ldots,L_t)$ is ``easier'' than
that of $\dim T/I(S_n(d))$.   We explored a number of ways
one could pick the $L_i$'s (e.g., picking forms at random, making use of the
symmetry), but no method allowed us to improve existing bounds, even with
our extensive computer resources.

\noindent
{\bf Acknowledgements.}
We thank the anonymous referees for their useful suggestions.
This work was made possible by the facilities of the Shared Hierarchical 
Academic Research Computing Network (SHARCNET:www.sharcnet.ca) and 
Compute/Calcul Canada.   
The first author was supported by
an NSERC USRA and the second author by an NSERC Discovery Grant.



\end{document}